\documentstyle[12pt]{article}

\newtheorem{lem}{Lemma}[section]
\newtheorem{rem}[lem]{Remark}
\newtheorem{thm}[lem]{Theorem}

\newtheorem{cj}[lem]{Conjecture}

\newcommand{\be}{\begin{equation}}
\newcommand{\ee}{\end{equation}}
\begin{document}
\setlength{\unitlength}{1mm}
\baselineskip 8mm
\thispagestyle{empty}
\title{A Generalization of Tepper's Identity}
\author
{Mortaza Bayat and Hossein Teimoori Faal
\vspace{.25in}\\
Department of Mathematics,
Zanjan Branch,\\ Islamic Azad University,
Zanjan, Iran  
\vspace{.25in}\\
Department of Mathematics and Computer Science,\\
Allameh Tabataba'i University, Tehran, Iran
}

\date{6 May 2021}

\maketitle
\begin{abstract}

In this paper, we first give 
a simple combinatorial proof of Tepper's identity. Then, as a by product of this interesting identity we present another 
proof of the well-knonw Wilson's identity in number theory. 
Finally, we obtain a generalization of Tepper's identity for any polynomial with real coefficients.    
\end{abstract}

\newpage

\section{Introduction}

There is no doubt that Pascal's triangle is one of the most beautiful triangular numerical array in mathematics. 
One can explore many algebraic, geometric and number-theoretic patterns inside this well-known numeric array. 
\\
The Newton's binomial identity is one of the famous algebraic identity thta one may encounter. That is, for example, for the third row of Pascal's 
triangle, we have 
\[
1 x^{0} + 3 x^{1} + 3 x^{2} + 1 x^{3} = (1+x)^{3}. 
\] 

By a more careful inspection, we may notice that, for example, for the third row, we can see 
\[
1(x-0)^{3} -3(x-1)^{3} + 3(x-2)^{3} -1 (x-3)^{3} = 3!. 
\]
The above observation may lead us to the important question whether the above identity is a coincidence 
or it is just a special case of a general identity. 
Indeed, one can conjecture that the \emph{inner product}
of the $n$th row of \emph{Pascal triangle} with signs alternating between $+$ and $-$ and the row-vector 

$$
((x-0)^{n}, (x-1)^{n}, \ldots, (x-n)^{n} ),
$$ 
is equal to $n!$ for any non-negative integer number $n$. 
In other words, we get the following identity 

\begin{equation}\label{Tepper1}
	\sum_{k=0}^{n}
	(-1)^{k} {n \choose k } (x -k)^{n} = n!.
\end{equation}

This identity is known as \emph{Tepper's identity}, as M. Tepper 
conjectured this result from a consideration of the numerical
data which he gave in $[1]$. 
In the same year, C. T. Long gave a proof of formula (\ref{Tepper1})
in $[2]$. Later, F. J. Papp gave another proof of 
it by mathematical induction in $[3]$. This result also 
\emph{implicitly} derived in problem $20$ by Feller $[4]$. 
\\
The paper organization is as follows. We first give a simple proof 
of Tepper's identity by a \emph{combinatorial argument}. Next, 
by using this beautiful identity, we present another proof of 
\emph{Wilson's theorem}. Finally, we generalize 
Tepper's identity for any \emph{polynomial} with \emph{real} coefficients. 

\section{A Combinatorial Method}
The following interesting combinatorial problem has an 
important role for the investigation of the \emph{Stirling} numbers 
of the \emph{first} kind. 
As the proof of Tepper's identity and it's generalization 
is closely related to this problem, we review 
a \emph{variant} of the original problem and it's solution here $[5]$.
\\
{\bf Problem (a).} 
Let a train have $n$ wagons. Now, if we \emph{randomly} choose 
a wagon, compute the number of ways in which \emph{exactly}
$r$ wagons will be occupied. 
\\
{\bf Problem (b).} Using the solution of the above problem
compute the following \emph{summation}:
\begin{equation}
	{n \choose 1} 1^{p}
	- {n \choose 2} 2^{p} +
	{n \choose 3} 3^{p} \cdots + 
	(-1)^{n-1}{n \choose n} n^{p
	}, \hspace{1cm} ( 1 \leq p  \leq n).
\end{equation} 

\begin{rem}
	The
	above \emph{statements} 
	are \emph{equivalent} to the following 
	interesting problem in \emph{physics}. 
	A \emph{sensor}
	contains $n$ \emph{receivers} and 
	receives a \emph{flux} of 
	$p$ particles. If the \emph{probability} of 
	receiving particles for any receiver is the same, then
	compute the probability that these particles will
	hit exactly $r$ receivers. 		  
\end{rem}

{\bf Solution (a).}
Suppose $A_{i}$
is the number of ways 
in which 
the $i$th wagon can be empty ($1 \leq i \leq n$).  
Then the number undesirable ways are 
$\vert A_{1} \cup A_{2} \cup \cdots \cup A_{r}  \vert$.
Hence, by the \emph{inclusion - exclusion principle}, we obtain
\begin{eqnarray}
	\vert A_{1} \cup A_{2} \cup \cdots \cup A_{r}  \vert
	& = &
	\sum_{i} \vert A_{i} \vert - 
	\sum_{i < j} \vert A_{i} \cap  A_{j} \vert + 
	\sum_{i < j < k} \vert A_{i} \cap  A_{j} \cap A_{k} \vert
	\nonumber\\
	& + &
	\cdots + (-1) ^ {n-1}
	\vert A_{1} \cap  A_{2} \cdots \cap A_{r} \vert
	\nonumber\\
	& = & 
	\sum_{i} (r-1)^{p} -  \sum_{i < j} (r-2)^{p}
	+ \cdots + (-1) ^{n-1} (r-r)^{p}
	\nonumber\\
	& = &
	{r \choose 1} (r-1)^{p}
	- {r \choose 2} (r-2)^{p} +
	\cdots + 
	(-1)^{r-2}{r \choose r-1} 1^{p
	}. 
	\nonumber
\end{eqnarray}
Since the number of all \emph{possible} cases is 
$r^{p}$, then the number of desirable ways is 
equal to
$$
r^{p} - 
{r \choose 1} (r-1)^{p}
+ {r \choose 2} (r-2)^{p} -
\cdots + 
(-1)^{r-1}{r \choose r-1} 1^{p
}. 
$$
Finally, since the number of ways of choosing
$r$ wagons is ${n \choose r}$, then the \emph{total} 
number of the desirable ways is 
\begin{equation}\label{equt2}
	{n \choose r} \left( 
	r^{p} - 
	{r \choose 1} (r-1)^{p}
	+ {r \choose 2} (r-2)^{p} -
	\cdots + 
	(-1)^{r-1}{r \choose r-1} 1^{p
	}
	\right).
\end{equation} 
{\bf Solution (b).}
Let $r=n$ in 
$(\ref{equt2})$. If $p < n$, the number of desirable ways
computed in part $(a)$ is equal to \emph{zero} and considering 
the well-known identity 
${n \choose r} = {n \choose n-r}$, we obtain
\begin{equation}\label{equt3}
	{n \choose 1} (1)^{p}
	- {n \choose 2} (2)^{p} +
	\cdots + 
	(-1)^{n-1}{n \choose n} n^{p
	}= 0.
\end{equation}
Now if we choose $n=p$ in (\ref{equt2}), 
then for \emph{every} passenger
we have \emph{exactly} one wagon and  
consequently the number of desirable 
ways is $n!$. Finally, considering 
$n=p$ and $r=p$ in (\ref{equt2}) and after simplifications, we have
\begin{equation}\label{equt4}
	{p \choose 1} (1)^{p}
	- {p \choose 2} (2)^{p} +
	\cdots + 
	(-1)^{p-1}{p \choose p} p^{p
	}= (-1)^{p-1} p!.
\end{equation}

\section{A Simple Proof of Tepper's Identity}
In this section, we give a simple proof of 
Tepper's identity, using the results of the previous 
section. 
\\
By substituting the \emph{Newton binomial expansion}
for $(x-k)^{n}$ in the left-hand side of 
(\ref{Tepper1}), we have
\begin{equation}\label{equt5}
	\sum_{k=0}^{n}
	(-1)^{k} {n \choose k } (x -k)^{n} = \sum_{j=0}^{n}
	x^{j} 
	\left( 
	\sum_{k=0}^{n}
	{n \choose k } (-k)^{n-j}
	\right). 
\end{equation}

Now, relations (\ref{equt3}) and 
(\ref{equt4}) lead us to 
\begin{eqnarray}
	\sum_{k=0}^{n}(-1)^{k}
	{n \choose k } (-k)^{n-j} & =  & 0, \hspace{1cm} (1 \leq j \leq n).\\
	\sum_{k=0}^{n}
	(-1)^{k} {n \choose k } (-k)^{n}& = & n!.
\end{eqnarray}
Thus, by substituting the above relations 
in (\ref{equt5}), we obtain Tepper's identity. 
As an immediate consequence of Tepper's identity, 
we obtain another proof of the \emph{Wilson's theorem}.

\begin{thm}[Wilson's theorem]
	Let 
	$p$ be an odd prime number. Then, we have
	$$
	(p-1)! \equiv -1\hspace{0.5cm} (mod~ p).
	$$
\end{thm}

\begin{pf}
	
Put $x=0$ and $n=p-1$ in Tepper's identity (\ref{Tepper1}). Then, we get
\begin{equation}
(p-1)! =
\sum_{k=0}^{p-1}
(-1)^{k}{p-1 \choose k} k^{p-1}.\nonumber 
\end{equation}

Now, using \emph{Fermat's little theorem}, we obtain
\begin{equation}
(p-1)! \equiv 
\sum_{k=1}^{p-1}
		(-1)^{k}{p-1 \choose k} \hspace{0.5cm} (mod~ p),
\end{equation}
or equivalently, 
\begin{equation}
		(p-1)! \equiv 
		-{p-1 \choose 0} +
		\sum_{k=0}^{p-1}
		(-1)^{k}{p-1 \choose k} \hspace{0.5cm} (mod~ p),\nonumber
\end{equation}
	which leads us to
	$$
	(p-1)! \equiv -1 \hspace{0.5cm} (mod~ p),
	$$
	as required. 
\end{pf}

\section{A Generalization of Tepper's Identity}

As we already saw in the previous section, the Newton's binomial 
expansion of $(x-k)^{n}$ had a \emph{key role} in
proving Tepper's identity. Now, we show that the \emph{identity}
is also valid for any polynomial of degree $n$.  
To do this, 
we first prove the following lemma by the same argument stated in 
$[6]$. 

\begin{lem}
	For any 
	real number $x$, and natural numbers 
	$m$ and $n$ such that, $0 \leq m < n$, we have
	\begin{equation}\label{lem1}
		\sum_{k=0}^{n}
		(-1)^{k} {n \choose k } (x -k)^{m} = 0.
	\end{equation}
\end{lem}

\begin{pf}
	By considering 
	the coefficient $x^{j}$, in the binomial expansion 
	of $(x-k)^{m}$ for $0 \leq j \leq m$, on the 
	left-hand side of (\ref{lem1}) and relation (\ref{equt3}) we get
	the desired result.
\end{pf}

\begin{rem}
	The 
	above lemma is also true for any polynomial $P(x)$ of degree 
	$m$, provided that $m < n$; that is,
	\begin{equation}
		\sum_{k=0}^{n}
		(-1)^{k} {n \choose k } P(x -k) = 0.
	\end{equation} 
\end{rem}

\begin{thm}[The Generalized Tepper's Identity]
Suppose $n$ is a natural number and P(x) is any 
	polynomial of degree $n$, as follows

$$
P(x) = a_{n} x^{n}
+ a_{n-1} x^{n-1} +
\cdots + a_{1} x^{1} + a_{0}, 
$$
in which $a_{i} \in \mathbf{R}$ for $(0 \leq i \leq n)$. 
Then, we have

\begin{equation}
		\sum_{k=0}^{n}
		(-1)^{k} {n \choose k } P(x -k) = a_{n} n!.
\end{equation}

\end{thm}

\begin{pf}
	By considering 
	Tepper's identity and lemma (\ref{lem1}), 
	we have
	\begin{eqnarray}
		\sum_{k=0}^{n}
		(-1)^{k} {n \choose k } P(x -k)
		& = & 
		\sum_{i=0}^{n} 
		\left( 
		\sum_{k=0}^{n} (-1)^{k} {n \choose k } (x-k)^{i}
		\right) 
		\nonumber \\
		& = & a_{n} 
		\left( 
		\sum_{k=0}^{n} (-1)^{k} {n \choose k } (x-k)^{i}
		\right) 
		\nonumber\\
		& = & 
		a_{n} n! \nonumber.
	\end{eqnarray}
\end{pf}

\begin{rem}
	Another 
	interesting
	proof
	of 
	Tepper's identity and 
	lemma (\ref{lem1}) is given 
	in $[7]$, using 
	Pascal's functional matrix. Also, several interesting 
	combinatorial identities are given in $[8]$.
We also come up with the following conjecture, as a generalization of 
\emph{Tepper's identity}

\begin{cj}
For any positive integer $l$ and any real polynomial $P(x)$ of 
degree $n$ with leading coefficient $a_{n}$, we have  

\begin{equation}
\sum_{k=0}^{n}
(-1)^{k} {n \choose k } P(x -lk) = a_{n} l^{n} n!.
\end{equation}
	
\end{cj}

\end{rem}



\end{document}